\documentstyle[12pt]{article}

\newcommand{\const}{\mathop{\rm const}\limits}

\newcommand{\Var}{\mathop{\rm Var}\limits}

\begin{document}

\begin{center}

{\bf  SUFFICIENT CONDITIONS FOR EXPONENTIAL TIGHTNESS } \\

\vspace{4mm}

{\bf  OF NORMED SUMS OF INDEPENDENT RANDOM FIELDS } \\

\vspace{4mm}

{\bf  IN THE SPACE OF CONTINUOUS FUNCTIONS } \\

\vspace{4mm}

 $ {\bf E.Ostrovsky^a, \ \ L.Sirota^b} $ \\

\vspace{4mm}

$ ^a $ Corresponding Author. Department of Mathematics and Computer science, \\
Bar-Ilan University, 84105, Ramat Gan, Israel.\\

\vspace{4mm}

E - mail: \ galo@list.ru \  eugostrovsky@list.ru\\

\vspace{4mm}

$ ^b $  Department of Mathematics and computer science. Bar-Ilan University, 84105,\\ Ramat Gan, Israel.\\

\vspace{4mm}

E - mail: \ sirota3@bezeqint.net \\

\vspace{5mm}
                    {\bf Abstract.}\\

 \end{center}

 \vspace{4mm}

  We formulate and prove in this report some sufficient conditions for {\it exponential tightness (ET)} of a family
 of independent identical distributed (i.i.d.) random fields (r.f.) (processes) in the space of continuous functions defined on
 certain metric compact set.\par
  We will use the entropic method as well as too modern approach through the majorizing measures (generic chaining).\par

 \vspace{4mm}


\vspace{4mm}


\vspace{4mm}

\section{Introduction. Statement of the problem. Notations. Definitions. }

\vspace{3mm}

 Let $ (\Omega, F, {\bf P } ) $ be a probability space,  $  (T = \{t\},d) $ be compact metric space
relative the semi - distance function $ d = d(t,s), \ t,s \in T, \  C(T,d) $ be the Banach
space of all continuous functions:  $ f:  \ T \to R  $  with ordinary uniform norm

$$
||f|| = \max_{t \in T} |f(t)|,
$$

 $ \xi = \xi(t), \ t \in T $ be separable random process (field),    $  \xi_i(t), \ i = 1,2, \ldots $ be
independent copies of $  \xi(t),  $

$$
 S_n(t) = n^{-1}\sum_{i=1}^n \xi_i(t).
$$

\vspace{3mm}

{\bf  Definition 1.1. Exponential Tightness. } We will say as usually that the family $ \{S_n(t) \}  $ satisfies the exponential tightness condition
(ETC) in the Banach space $  C(T,d) $, briefly: $  \{ S_n(\cdot) \} \in ETC, $
 if the r.f. $ \xi(\cdot)  $ and following all the r.f. $ S_n(\cdot) $ are $ d \ - $ continuous  with probability one,  and
 if for each $ v = \const > 0 $ there exists a compact set $  K = K(v) $ in the space $ C(T,d) $  such that

$$
\lim \sup_{n \to \infty} \frac{1}{n} \log {\bf P} (S_n(\cdot) \notin K(v) ) \le -v. \eqno(1.1)
$$

 This notion play a very important role in the theory of Large Deviations (LD) in the Banach space, see
\cite{Acosta1} - \cite{Acosta3},  \cite{Bahadur1}, \cite{Borovkov1}, \cite{Dembo1}, \cite{Deuschel1},
\cite{Freidlin1}, \cite{Piterbarg1}, \cite{Puhalskii2}, \cite{Schied2}, \cite{Varadhan1} and following in Statistics
\cite{Bahadur1}, \cite{Bryc1}, \cite{Chaganty1}, \cite{Choirat1}, \cite{Lindsay1}, \cite{Ostrovsky103}, \cite{Puhalskii1};
in the theory of Random Processes \cite{Piterbarg1}, \cite{Acosta3}; in Insurance and Finance \cite{Kluppelberg1};
and still in Philology \cite{Ostrovsky103}.\par

\vspace{3mm}

{\bf Our purpose in this report is to formulate and prove some sufficient conditions for exponential tightness in the
Banach space of continuous functions. }

\vspace{3mm}

 Some previous results in this direction may be found in the articles
\cite{Acosta1} - \cite{Acosta2}, \cite{Ney1}, \cite{Piterbarg1}, \cite{Puhalskii3}, \cite{Schied2}, \cite{Schied1}, \cite{Schied3}.\par
 Note that it is written in the name of the article \cite{Schied2} "Criteria for exponential tightness in path spaces", but
 in fact it is formulated and proved only sufficient conditions for exponential tightness in path space. \par

\vspace{3mm}

{\it Roughly speaking, we will prove that "almost everywhere", i.e. under simple natural additional conditions,
  when  the sequence of r.f. $ n^{1/2} S_n(\cdot) $ satisfies the Central
Limit Theorem (CLT) in the space $  C(T,d), $  then it satisfies also the Exponential Tightness (ET) in this space.} \par

\vspace{3mm}

\section{ Entropy conditions for exponential tightness.}

\vspace{3mm}

 We can and will suppose without loss of generality that the r.f. $ \xi = \xi(t) $ is non-zero, centered: $ {\bf E } \xi(t) = 0, \ t \in T. $\par
We suppose  also that  $ \xi(t)  $ satisfies the uniformly in $  t \in T  $  the following strengthening of Cramer's  condition:

$$
\forall  \mu = \const > 0  \ \exists x_0 > 0 \ \forall x > x_0 \Rightarrow \sup_{t \in T} {\bf P} (|\xi(t)| >x) \le \exp(-\mu x). \eqno(2.1)
$$

Introduce the following important function:

$$
\phi(\lambda)  \stackrel{def}{=} \sup_{t \in T} \max_{\pm} \log {\bf E} \exp(\pm \lambda \xi(t)). \eqno(2.2)
$$
 It follows from condition (2.1) that the function $ \phi(\lambda) $ is finite for all the values $ \lambda, \ \lambda \in R,  $ and
is even, is logarithmical convex  and such that

$$
\phi(0) = 0, \ \phi'(0) = 0,  \ |\lambda| \le 1 \ \Rightarrow C_1 \lambda^2 \le \phi(\lambda) \le C_2 \lambda^2 \eqno(2.3)
$$
for some finite positive constants $ C_1, C_2. $ \par
 Define on the basis of the function $ \phi(\cdot) $ the next function

$$
\chi(\lambda) \stackrel{def}{=} \sup_n \ [n \cdot \phi(\lambda/\sqrt{n})]. \eqno(2.4)
$$
 This function obeys at the same properties as the source function $ \phi(\cdot); $ for instance, the finiteness of the
function $ \chi(\lambda) $  for all the values $ \lambda \in R $ follows immediately from the finiteness of the function
$ \phi $ and from the relations (2.3). \par

 The function $ \chi = \chi(\lambda) $ generated the so-called Banach space $ B(\chi) = B(\chi; \Omega) $ as follows.
 We say that the {\it centered} random variable (r.v) $ \eta  $
belongs to the space $ B(\chi), $ if there exists some non-negative constant
$ \tau \ge 0 $ such that

$$
\forall \lambda \in R  \ \Rightarrow
{\bf E} \exp(\lambda \eta) \le \exp[ \chi(\lambda \ \tau) ]. \eqno(2.5).
$$
 The minimal value $ \tau $ satisfying (2.5) is called a $ B(\chi) \ $ norm
of the variable $ \xi, $ write
 $$
 ||\eta||B(\chi) = \inf \{ \tau, \ \tau > 0: \ \forall \lambda \ \Rightarrow
 {\bf E}\exp(\lambda \eta) \le \exp(\chi(\lambda \ \tau)) \}. \eqno(2.6)
 $$
 These spaces are Banach spaces and are very convenient for the investigation of the r.v. having an
exponential decreasing tail of distribution, for instance, for investigation of the limit theorem,
the exponential bounds of distribution for sums of random variables,
non-asymptotical properties, problem of continuous of random fields,
study of Central Limit Theorem in the Banach space etc.\par
 The detail investigation of these spaces see in \cite{Kozatchenko1},  \cite{Ostrovsky101}, chapter 1.\par

\vspace{3mm}

 Note that both the functions $ \phi, \ \chi, $ as well as the norm $ || \cdot||B(\chi)  $ are {\it natural}, i.e.
generated only by the values of $ \xi(t), \ t \in T. $  We can define by means of these functions and norm
define the following variable $ \sigma := \sup_{t \in T} ||\xi(t)||B(\chi) $ and distance (more exactly, semi-distance):

$$
d_{\phi}(t,s) := ||\xi(t) - \xi(s)||B(\phi), \eqno(2.7)
$$
which is also natural. \par

 Notice that the norms $ ||\cdot||B(\phi) $ and $ ||\cdot||B(\chi)   $  are in general case not equivalent, but if
 $ \sup_{t \in T} ||\xi(t)||B(\phi) < \infty,  $ then $ \sigma < \infty. $\par

 For example, if the r.f. $ \xi(t)  $ is centered and Gaussian with $ \max_{t \in T} \Var(\xi(t)) = 1, $ then

$$
\sigma = 1, \ \hspace{7mm} \phi(\lambda) = \chi(\lambda) = \lambda^2/2
$$
and correspondingly

$$
d_{\chi}(t,s)  = \left[ \Var(\xi(t) - \xi(s))   \right]^{1/2}. \eqno(2.8)
$$

 The distance  $ d_{\chi}(t,s) $ in  (2.8) is called Dudley, or Dudley - Fernique distance, see \cite{Dudley1}, chapter 2,3;
\cite{Fernique1}.\par

 We denote as ordinary by $ N(T, d_{\chi}, \epsilon) $ the minimal amount of  $ d_{\chi} \ - $ balls of a radii  $ \epsilon > 0 $
which  cover all the set $  T. $ Recall that it follows from Hausdorff's theorem that the set $  T  $ relative the distance
$ d_{\chi}  $ is compact iff it is bounded, closed and for arbitrary  $ \epsilon  > 0 \ \Rightarrow  N(T, d_{\chi}, \epsilon) < \infty. $ \par

 The function

 $$
 H(T, d_{\chi}, \epsilon)  =  \ln N(T, d_{\chi}, \epsilon)
 $$
is called an {\it  entropy} of the set $  T  $ relative the semi-distance $ d_{\chi}. $ \par

 Introduce also the following functions and variables.

$$
\chi^*(x) = \sup_{\lambda \in R} ( \lambda |x| - \chi(\lambda))
$$
the Young - Fenchel, or Legendre transform of the function $ \chi(\cdot); $

$$
Y(z) := \exp\left(\chi^*(z) \right) - 1
$$
be a Young - Orlicz function generated by $ \chi(\cdot);  $

$$
J:= \int_0^{\sigma}  Y^{-1} \left( N(T, d_{\chi}, \epsilon) \right) \ d \epsilon. \eqno(2.9)
$$

 The integral $  J  $ is called {\it entropy  } integral; the condition $ J < \infty  $ is said to be {\it  entropy condition}. \par

\vspace{4mm}

{\bf Theorem 2.1.}  {\it Assume $  J < \infty.  $
Then the sequence $ S_n(\cdot) $  is exponentially tight in the Banach space } $ C(T, d_{\phi}). $ \par

\vspace{3mm}

{\bf Proof.} \par

{\bf 0.} We will use hereafter the following famous result belonging to Alejandro de Acosta \cite{Acosta1} - \cite{Acosta2}
on sums of Banach space-valued random variables. Let $ E $ be a Banach space with a norm $ ||\cdot|| $ , and let
$ \eta_1, \eta_2, \ldots $ be independent and identically distributed $ E \ - $ valued random variables satisfying the condition

$$
\forall \lambda = \const \ge 0 \ \Rightarrow {\bf E} \exp(\lambda ||\eta_1 ||) < \infty.
$$

Define

$$
X_n = \frac{1}{n} \sum_{j=1}^n \eta_j.
$$

A. de Acosta proved  that $ \{  X_n  \} $ is exponentially tight.\par

{\bf 2.} We can and will suppose

$$
\lim_{\lambda \to \infty} \phi(\lambda)/\lambda = \lim_{\lambda \to \infty} \chi(\lambda)/\lambda = \infty,
$$
since the opposite case is trivial: for some finite constant $  C  $

$$
{\bf P} ( \sup_{t \in T} |\xi(t)| \le C  ) = 1.
$$

 Further, let us consider the random field, more exactly, the family of random fields
 $$
 \zeta_n(t) = n^{1/2} S_n(t) = n^{-1/2}\sum_{i=1}^n \xi_i(t).  \eqno(2.10)
 $$
We have:

$$
\sup_n \sup_{t \in T} ||\zeta_n(t)||B(\chi) = \sigma < \infty;  \eqno(2.11a)
$$

$$
\sup_n \sup_{t \ne s} \left| \left| \frac{ \zeta_n(t) - \zeta_n(s)}{d_{\phi}(t,s) } \right| \right|B(\chi) =:C_3 < \infty. \eqno(2.11b)
$$

  Since the entropy condition $  J < \infty $ is satisfied, there exists a new continuous relative the semi-distance $ d_{\phi} $
semi-distance $  \rho = \rho(t,s) $ such that

$$
{\bf P} \left\{  \sup_{t \ne s} \left| \frac{\zeta_n(t) - \zeta_n(s)}{\rho(t,s)} \right|  > u   \right\} \le
\exp \left( - \chi^*(C u)  \right),  \eqno(2.12)
$$
see  \cite{Dmitrovsky1} or \cite{Ostrovsky101}, chapter 4, section 4.3. \par

 Therefore

 $$
{\bf P} \left\{  \sup_{t \ne s} \left| \frac{S_n(t) - S_n(s)}{\rho(t,s)} \right|  > u   \right\} \le
\exp \left( - \chi^*(C u \sqrt{n})  \right).  \eqno(2.13)
$$

 This completes the proof of theorem 2.1. \par

\vspace{3mm}

\section{ Majorizing measure conditions for exponential tightness.}

\vspace{3mm}

Let $ m(\cdot) $ be probabilistic Borelian measure defined on the set $ T,$  and
let $ \Phi = \Phi(u), \ u \ge 0 $ be the strictly increasing continuous non-negative Young-Orlicz function such that

$$
\Phi(u) = 0 \ \Rightarrow u = 0; \ \forall \lambda > 0 \ \Rightarrow \lim_{u \to \infty} \frac{\Phi(u)}{\exp(\lambda u)} = \infty.
\eqno(3.1)
$$

 Let us define also the following important distance function: $ w(x_1, x_2) = $
 $$
  w(x_1, x_2; V) = w(x_1, x_2; V, m ) = w(x_1, x_2; V, m,\Phi) = w(x_1, x_2; V, m,\Phi,d) \stackrel{def}{=}
 $$

$$
 6 \int_0^{d(x_1, x_2)} \left\{ \Phi^{-1} \left[ \frac{4V}{m^2(B_d(r,x_1))} \right] +
\Phi^{-1}  \left[ \frac{4V}{m^2(B_d(r,x_2))} \right] \right\} \ dr, \eqno(3.2)
$$
 where $ B_d(x,r) $ is $  d \ - $ ball in the set $  X: $

$$
B_d(x,r) = \{ y, \ y \in X, \ d(x,y) \le r   \}, \ r \le D := \sup_{x,y} d(x,y).
$$

  The triangle inequality and other properties of the distance function $ w = w(x_1, x_2) $ are proved in
\cite{Kwapien1}.\par

\vspace{3mm}

{\bf Definition 3.1. } (See  \cite{Kwapien1}). The measure $ m  $ is said to be
{\it minorizing measure } relative the distance $ d = d(x_1,x_2), $ if for each values $ x_1, x_2 \in X, \ V \in (0, \infty) $
 $ \ w(x_1,x_2; V) < \infty. $\par

\vspace{3mm}

{\bf Definition 3.2. } (See  \cite{Kwapien1}). The measure $ m  $ is said to be
{\it majorizing measure } relative the distance $ d = d(x_1,x_2), $ if

$$
\sup_{x,y \in X}  w(x_1,x_2; D) < \infty.
$$

\vspace{3mm}

  We will denote the Orlicz norm constructed by means of the function $  \Phi $
 of a random variable $ \kappa $ defined on our probabilistic space $ (\Omega, {\bf P}) $ as $ |||\kappa|||L(\Omega,\Phi) $   or for simplicity
 $ |||\kappa|||\Phi. $  \par
  We introduce  the so-called {\it natural} distance $ \rho_{\Phi}(x_1,x_2) $ as follows:

  $$
  d_{\Phi} = d_{\Phi}(x_1,x_2):= |||\xi(x_1) - \xi(x_2)|||L(\Omega,\Phi), \ x_1,x_2 \in X.  \eqno(3.3)
  $$

\vspace{4mm}

{\bf Theorem 3.1.} {\it  Suppose  for some described above Young - Orlicz function} $ \Phi(\cdot) $

$$
\inf_{t \in T} |||\xi(t)|||L(\Omega,\Phi) < \infty, \eqno(3.4)
$$
{\it and that there exists  the probabilistic Borelian majorizing measure $ m(\cdot) $   on the set $  X  $
relative the distance } $  d_{\Phi}(\cdot, \cdot). $ \par

{\it  Then the sequence $ S_n(\cdot) $  is exponentially tight in the Banach space } $ C(T, d_{\Phi}). $ \par

\vspace{4mm}

{\bf Proof.} Fix arbitrary point $ t_0 \in T $ for which  $ ||\xi(t_0)|||L(\Omega,\Phi) < \infty. $ Further,
there exists a non-negative random variable $ \theta $ belonging to the space  $  L(\Omega, \Phi) $ and continuous
relative the distance $ d_{\Phi} $ (and bounded) a new non-random distance $ w = w(t_1, t_2) $ for which

$$
|\xi(t_1) - \xi(t_2)| \le \theta \cdot w(t_1, t_2), \eqno(3.5)
$$
\cite{Ostrovsky106}; see  also \cite{Ostrovsky104}, \cite{Ostrovsky105}.  Following, the r.f.
$  \xi = \xi(t), \ t \in T  $ is $ d_{\Phi} \ - $ continuous with probability one. \par

\vspace{3mm}

 Further, it follows from the expression

$$
\xi(t) = \xi(t_0) + (\xi(t) - \xi(t_0)),
$$
and from the inequality (3.5), properties of the point $ t_0 $ that

$$
\sup_{t \in T} |\xi(t)| \in L(\Omega, \Phi).
$$

 We conclude taking  into account the properties of the function $ \Phi  $ (3.1)  that

$$
\forall \lambda \ge 0 \ \Rightarrow \ {\bf E}\exp \left( \lambda \sup_{t \in T} |\xi(t)|  \right) < \infty. \eqno(3.6)
$$
 The proposition of our theorem (3.1) follows now from one of the main results of an articles Alejandro de Acosta \cite{Acosta1}
  - \cite{Acosta2}. \par

\vspace{3mm}

\section{ Exponential tightness in  H\"older spaces.}

\vspace{3mm}

 Let $  \omega = \omega(t,s) $ be another continuous non-random semi-distance on the set $  T. $ The H\"older's  space
 $ H(T,\omega) = H(\omega) $ is defined  as a set of all functions $f, \ f: T \in R $ with finite norm

 $$
 ||f||H(\omega) \stackrel{def}{=} |f(t_0)| + \sup_{t,s: t \ne s } \left[\frac{|f(t) - f(s)|}{\omega(t,s)}  \right], \eqno(4.1)
 $$
where $  t_0 $ is certain fixed non-random point of the set $ T. $ \par

 We derive in this section some sufficient conditions for Exponential Tightness in H\"older spaces. \par

\vspace{4mm}

{\bf Theorem 4.1.} {\it We retain all the notations and conditions of theorem 2.1, in particular, the condition $  J < \infty. $  \par
Let $ \rho = \rho(t,s) $ be a continuous relative $ d_{\Phi} $  semi-distance  defined in the second section. \par
 The sequence $ S_n(\cdot) $  is exponentially tight in the Banach space } $ H(T, \rho). $  \par

\vspace{3mm}

{\bf Proof} is very simple. It follows from relations (2.12) and (2.13) that

$$
\forall \lambda \ge 0 \ \Rightarrow \sup_n {\bf E} \exp\left(\lambda ||\zeta_n(\cdot)||H(\rho)  \right) < \infty. \eqno(4.2)
$$
 It remains to refer the famous cited result of A.de Acosta \cite{Acosta1}.\par
Analogously may pe grounded the following proposition. \par

\vspace{4mm}

{\bf Theorem 4.2.} {\it We retain all the notations and conditions of theorem 3.1.  Let  } $ w = w(t,s) $
{\it be the semi-distance defined in the third section. We claim that the sequence $ S_n(\cdot) $  is exponentially tight
in the Banach space } $ H(T, w). $   \par

\vspace{3mm}

 Another exponential estimates for distribution of maximum for continuous random fields may be found in the article
\cite{Ostrovsky104}.\\


\vspace{4mm}

\begin{thebibliography}{99}

\vspace{3mm}

 \bibitem{Acosta1}
 {\sc de Acosta, Alejandro.} {\it Large deviations for vector-valued Levy processes.} Stochastic Process.
 Appl. 51 (1994), 75 \ - \ 115.

\bibitem{Acosta2}
{\sc de Acosta, Alejandro.}  {\it Exponential Tightness and Projective Systems in Large Deviation Theory.}
Festschrift for Lucien Le Cam. Springer, New York. 143 \ - \ 156, 1997.

 \bibitem{Acosta2}
 {\sc de Acosta, Alejandro.} {\it Large deviations for vector-valued Levy processes.} Stochastic Process.
 Appl. 51 (1994), 75 \ - \ 115.


\bibitem{Acosta3}
{\sc de Acosta, Alejandro.}   {\it A general non-convex large deviation result with applications to stochastic
equations.} Probab. Theory Relat. Fields, {\bf 118,} (2000), 483 \ - \ 521.

\bibitem{Arnold1}
{\sc Arnold L. and Imkeller P.} {\it On the spatial asymptotic Behavior of stochastic
Flows in Euclidean Space. } Stoch. Processes Appl., 62(1), (1996), 19-54.


\bibitem{Bahadur1}
{\sc Bahadur, R., Zabell, S., and Gupta, J.} (1980). {\it Large deviations, tests, and estimates.}
 In I.M. Chaterabarli (ed.), Asymptotic Theory of Statistical Tests and Estimation, pp. 33 \ - \ 64. New York: Academic Press.
Mathematical Reviews (MathSciNet): MR571334

\bibitem{Borovkov1}
{\sc Borovkov A.A. and Mogul'skii A.A.} {\it  Probabilities of large deviations in topological spaces, I, II.  }
Sibirsk. Math. Zh., {\bf 19,} (1978), 988 \ - \ 1004; {\bf 21 (5),} (1980), 12 \ - \ 26;  English transl. in Siberian  Math. J.
{\bf 19,} (1978), {\bf 21,} (1980).

\bibitem{Bryc1}
{\sc Bryc W. } (1990). {\it Large deviations by the asymptotic value method. }
In M.A. Pinsky (ed.), Diffusion Processes and Related Problems in Analysis, Vol. 1, pp. 447 \ - \ 472. Birkhäuser.

\bibitem{Chaganty1}
{\sc Chaganty, N.} (1993). {\it Large deviations for joint distributions and statistical applications.}
Technical Report TR93-2, Department of Mathematics and Statistics, Old Dominion University, Norfolk, VA.


\bibitem{Choirat1}
{\sc Choirat  Ch. and  Seri R.}
{\it Estimation in Discrete Parameter Models.}
Statistical Science, (2012,) Vol. 27, No. 2, $ 278 \ – \ 293. $
DOI: 10.1214/11-STS371.\\
 Also: arXiv:1207.5653v1 [stat.ME] 24 Jul 2012.


\bibitem{Deuschel1}
{\sc Deuschel, J. and Stroock, D.} (1989). {\it Large Deviations.}  Boston: Academic Press.
Mathematical Reviews (MathSciNet): MR997938

\bibitem{Dembo1}
{\sc Dembo A. and Zeitouni O.} {\it  Large deviations techniques and applications.} Jones and Bartlett, Boston, (1993).


\bibitem{Dmitrovsky1}
{\sc Dmitrovsky V.A.} {\it  On the  maximum distribution and local properties of paths  pre - Gaussian  random fields. }
Probab. Theory and Math. Statist., Kiev, KSU, (1981), {\bf 25}, 154 \ - \ 164.

\bibitem{Freidlin1}
{\sc Freidlin, M.I. and Wentzell, A. D.}
(1979) {\it Random Perturbations of Dynamical Systems.} Moscow: Nauka [in Russian]. English translation: Springer (1984).
Mathematical Reviews (MathSciNet): MR722136

\bibitem{Dudley1}
{\sc Dudley R.M.} {\it Uniform Central Limit Theorem.}  Cambridge University Press, 1999.

\bibitem{Fernique1}
 {\sc Fernique X.} (1975). {\it Regularite des trajectoires des
    function aleatiores gaussiennes.}  Ecole de Probablite de
    Saint-Flour, IV - 1974, Lecture Notes in Mathematic. {\bf 480}, 1 - 96,
    Springer Verlag, Berlin.

\bibitem{Kluppelberg1}
{\sc Kl\"uppelberg  C., Mikosch T.} {\it Large Deviations on Heavy-Tailed Random Sums  with applications in Insurance and Finance.}
 J. Appl. Probab., {\bf 34,} (1997),  293 \ - \ 308.

\bibitem{Kolmogorov1}
{\sc Kolmogorov, A. N. and Tikhomirov, V. M.}  (1959), {\it $ \epsilon $-entropy and $ \epsilon $ - capacity
of sets in a functional space. } Uspekhi Mat. Nauk, 14, 3; {\bf 86}.

\bibitem{Kozatchenko1}
{\sc Yu.V. Kozatchenko and E.I. Ostrovsky. } {\it Banach spaces of random
variables of subgaussian type.} Theory Probab. Math. Stat., Kiev,
1985,  42--56 (Russian).

\bibitem{Kwapien1}
{\sc Kwapien S. and Rosinsky J.} {\it Sample H\"older continuity of stochastic processes
and majorizing measures.} (2004). Seminar on Stochastic Analysis, Random
Fields and Applications IV, Progr. in Probab. 58, 155163. Birkh\"ouser, Basel.

\bibitem{Li1}
{\sc Li Y. } (2003). {\it A martingale inequality and large deviations.} Statist. Probab. Lett.
{\bf 62,} 317 \ - \ 321.


\bibitem{Lindsay1}
{\sc Lindsay, B. G. and Roeder, K.} (1987). {\it A unified treatment of integer parameter models.} J. Amer. Statist.
Assoc., {\bf 82}, 758 \ – \ 764. MR0909980

\bibitem{Meeden1}
{\sc Meeden, G. and Ghosh, M.} (1981).  {\it Admissibility in finite problems. } Ann. Statist., {\bf 9,} 846 \ – \ 852. MR0619287

\bibitem{Ney1}
{\sc Ney P. } (1984). {\it Convexity and large deviations. } Ann.
Probab., {\bf 12,}  903 \ – \ 906. MR0744245

\bibitem{Ostrovsky101}
{\sc  Ostrovsky E.I.} {\it Exponential Estimations for Random Fields.}
Moscow - Obninsk, OINPE, 1999 (in Russian).

\bibitem{Ostrovsky102}
{\sc Ostrovsky E.I.} (2002). {\it Exact exponential estimations for random field maximum
distribution.} Theory Probab. Appl. 45 v.3, 281 - 286.

\bibitem{Ostrovsky103}
{\sc Ostrovsky E., Sirota L., Zeldin A.}
{\it Rate of convergence in the maximum likelihood estimation for partial discrete parameter,
with applications to the cluster analysis and philology.}
arXiv:1402.6409v1 [math.ST] 26 Feb 2014

\bibitem{Ostrovsky104}
{\sc Ostrovsky E., Sirota L.} {\it Simplification of the majorizing measure method,
with development.} arXiv:1302.3202v1 [math.PR] 13 Feb 2013

\bibitem{Ostrovsky105}
{\sc Ostrovsky E., Rogover E.} {\it Exact exponential Bounds for the random Field
Maximum Distribution via the Majorizing Measures (Generic Chaining.)}\\
 arXiv:0802.0349v1 [math.PR] 4 Feb 2008

\bibitem{Ostrovsky106}
{\sc Ostrovsky E., Sirota L.}  {\it Uniform measures on the arbitrary compact metric spaces, with applications. }
arXiv:1403.5725v1 [math.FA] 23 Mar 2014

\bibitem{Piterbarg1}
{\sc Piterbarg V.I.  and Fatalov V.R.} {\it  The Laplace method for probability measures in Banach spaces.}
Russian Math. Surveys, 1995, {\bf  50,} 1151 \ - \ 1239.

\bibitem{Puhalskii1}
{\sc Puhalskii A. and  Spokoiny V.} {\it On large-deviation efficiency in statistical inference.}
 Bernoulli, Volume 4, Number 2 (1998), 203-272.

\bibitem{Puhalskii2}
{\sc Puhalskii, A. } (1991) {\it On functional principle of large deviations.} In V. Sazonov and T. Shervashidze (eds.),
New Trends in Probability and Statistics, Vol. 1, pp. 198-218. Utrecht: VSP/Moks'las.
Mathematical Reviews (MathSciNet): MR1200917

\bibitem{Puhalskii3}
{\sc Puhalskii, A.} (1993)  {\it On the theory of large deviations.} Theory Probab. Appl., 38(3), 490-497.
Mathematical Reviews (MathSciNet): MR1404664

\bibitem{Schied1}
{\sc Schied A. } {\it Sample path large deviations for super-Brownian motion.} Preprint, 1994.

\bibitem{Schied2}
{\sc Schied A.} {\it Criteria for exponential tightness in path spaces.}  Internet preprint,
1995.

\bibitem{Schied3}
{\sc Schied A.} {\it $ Gro \beta e $ Abweichungen f\"ur die Pfade der Super-Brownschen Bewegung.}
Bonner mathematische Schriften, Vol. 277 (1995).

\bibitem{Talagrand1}
{\sc Talagrand M.} (1996). {\it Majorizing measure: The generic chaining.}  Ann.
Probab., 24 1049 - 1103. MR1825156

\bibitem{Talagrand2}
{\sc Talagrand M.} (2005). {\it The Generic Chaining. Upper and Lower Bounds of
Stochastic Processes. } Springer, Berlin. MR2133757

\bibitem{Varadhan1}
{\sc Varadhan, S.R.S.} (1984) {\it Large Deviations and Applications.} Philadelphia: SIAM.
Mathematical Reviews (MathSciNet): MR758258


\bibitem{Varadhan2}
{\sc Varadhan, S.R.S.} (1966) {\it Asymptotic  probabilities  and differential equations.}
Comm. Pure Appl. Math.; {\bf 19,} (1966), 261 \ - 286.


\end{thebibliography}
\end{document}